\documentclass[
    ,final            
  ]
  {aipproc}

\layoutstyle{6x9}

\newtheorem{prop}{Proposition}
\newtheorem{theorem}{Theorem}
\newtheorem{lemma}{Lemma}
\newtheorem{coro}{Corollary}

\def\espaitemps{({\cal V},g)}
\def\varietat{{\cal V}}

\begin{document}

\title{Second-Order Symmetric Lorentzian Manifolds}

\classification{02.40.Ky, 04.20.-q, 04.50.+h}
\keywords      {Symmetric spaces, curvature invariants, parallel null 
vector fields, Mp-waves.}

\author{Jos\'e M. M. Senovilla}{
  address={F\'{\i}sica Te\'orica, Universidad del Pa\'{\i}s Vasco,
  Apartado 644, 48080 Bilbao, Spain.}
}

\begin{abstract}
Spacetimes with vanishing second covariant derivative of the Riemann tensor are studied. Their existence, classification and explicit local expression are considered. Related issues and open questions are briefly commented.
\end{abstract}

\maketitle


\section{Introduction}
Our aim is to characterize, as well as to give a full list of, the $n$-dimensional
manifolds $\varietat$ with a metric $g$ of Lorentzian signature such 
that the Riemann tensor $R^{\alpha}{}_{\beta\gamma\delta}$ of 
$\espaitemps$  {\em locally} satisfies the second-order condition 
\begin{equation}
\nabla_{\mu}\nabla_{\nu} R^{\alpha}{}_{\beta\gamma\delta}=0. 
\label{DDR}
\end{equation}
It is quite surprising that, hitherto, despite their simple 
definition, this type of Lorentzian manifolds have been hardly 
considered in the literature. Probably this is due to the classical results concerning these manifols in the proper Riemannian case, to the difficulties arising in other signatures, 
and to the little reward: only very special cases survive.

Apart from their obvious mathematical interest, from a physical point of view they are relevant in several respects: as a second local approximation to any spacetime (using for instance expansions in normal coordinates); as examples with a finite number of terms in Lagrangians; as interesting exact solutions for 
supergravity/superstring or M-theories; for invariant classifications; for solutions with parallel vector fields or spinors. 

A more complete treatment, with a full list of references, is given in 
\cite{2sym}.

\section{Symmetric spaces and its generalizations}
Semi-Riemannian manifolds satisfying 
(\ref{DDR}) are a direct generalization of the classical 
locally {\em symmetric} spaces which satisfy
\begin{equation}
\nabla_{\mu} R^{\alpha}{}_{\beta\gamma\delta}=0. \label{DR}
\end{equation}
These were introduced, studied and classified by E. Cartan 
\cite{C} in the proper Riemannian case\footnote{With a positive-definite metric.}, see e.g. \cite{C1,KN,H}, and later in 
\cite{CM,CW,CP} for the Lorentzian and general semi-Riemannian 
cases---see e.g. \cite{CP1,O} and references therein. They are themselves generalizations of the constant 
curvature spaces and, actually, there is a  
hierarchy of conditions, shown in Table \ref{ta}, that can be placed 
on the curvature tensor. In the table, the restrictions on the 
curvature tensor decrease towards the right and each class is 
{\em strictly} contained in the following ones. The table has been stopped 
at the level of semi-symmetric spaces, defined by the 
condition\footnote{(Square) round brackets enclosing indices indicate 
(anti-)symmetrization, respectively.} $\nabla_{[\mu}\nabla_{\nu]} R^{\alpha}{}_{\beta\gamma\delta}=0$, which were introduced also by Cartan \cite{C1} and studied in \cite{Sz,Sz1} as the natural generalization of symmetric 
spaces for the proper Riemannian case---see also \cite{BKV} and 
references therein.

\begin{table}
\begin{tabular}{|l|c|c|l|}
\hline
$R^{\alpha}{}_{\beta\gamma\nu}\propto 
\delta^{\alpha}_{\gamma}g_{\beta\nu}-\delta^{\alpha}_{\nu}g_{\beta\gamma}$ 
& $\nabla_{\mu} R^{\alpha}{}_{\beta\gamma\delta}=0$ &
$\nabla_{\mu}\nabla_{\nu} R^{\alpha}{}_{\beta\gamma\delta}=0$ & 
$\nabla_{[\mu}\nabla_{\nu]} 
R^{\alpha}{}_{\beta\gamma\delta}=0$ \\
\hline
constant curvature & symmetric & 2-symmetric & semisymmetric \\
\hline
\end{tabular}
\caption{The hierarchy of conditions on the Riemann tensor}
\label{ta}
\end{table}

Why was semisymmetry considered to 
be the natural generalization of local symmetry? And, why not going further on to higher derivatives of the Riemann tensor? The answer to both questions is actually the 
same: a classical theorem \cite{L,NO,Ta} states that in any proper 
Riemannian manifold
\begin{equation}
\nabla_{\mu_1}\dots \nabla_{\mu_k}R^{\alpha}{}_{\beta\gamma\delta}=0 
\hspace{1cm} \Longleftrightarrow  \hspace{1cm} \nabla_{\mu} 
R^{\alpha}{}_{\beta\gamma\delta}=0 \label{DkR}
\end{equation}
for any $k\geq 1$ so that, in particular, (\ref{DDR}) is strictly 
equivalent to (\ref{DR}) in proper Riemannian spaces. 
This may well be the reason why there seems to be no name 
for the condition (\ref{DDR}) in the literature. However, an 
analogous condition has certainly been used for the so-called 
$k$-recurrent spaces \cite{Tak,CS1};  thus, I will call the spaces satisfying (\ref{DDR}) {\em second-order symmetric}, or in short {\em 2-symmetric}, and more generally 
$k$-symmetric when the left condition in (\ref{DkR}) holds---see \cite{2sym} for further details.

\subsection{Results at \emph{generic} points}
As a matter of fact, the equivalence (\ref{DkR}) holds as well in 
``generic" cases of semi-Riemannian manifolds of any signature. For 
some results on this one can consult \cite{Ta,CS1}. By ``generic point" the following is meant: any point $p\in \varietat$ where the matrix 
$(R^{\alpha\beta}{}_{\gamma\delta})|_p$ of the Riemann tensor, 
considered as an endomorphism on the space of 2-forms $\Lambda_2(p)$, 
is non-singular. Then, for instance one can 
prove the following general result, see \cite{2sym} for a proof.
\begin{prop}
\label{2implies1}
For any tensor field $T$, and at \underline{{\em generic}} 
points, one has
$$
\stackrel{k}{\overbrace{\nabla\cdots\cdots\nabla}}\,  T =0 \Longleftrightarrow \nabla T =0
$$
for any $k\geq 1$.
\end{prop}
Of course, these results apply in particular to the Riemann tensor, and in fact sometimes even stronger results can be proven. For instance, one can prove a conjecture in \cite{CS1}, namely, that all $k$-symmetric (and also 
all $k$-recurrent) spaces are necessarily of constant curvature on a 
neighbourhood of any generic $p\in \varietat$. As a matter of fact, a slightly more general result is proven in \cite{2sym}:
\begin{theorem}
All {\em semi-symmetric} spaces are of constant curvature at generic 
points.
\end{theorem}
Therefore, there is little room for spaces (necessarily of 
non-Euclidean signature) which are $k$-symmetric but {\em not} 
symmetric nor of constant curvature. It is remarkable that there have 
been many studies on 2-recurrent spaces, but surprisingly enough the assumption that they are {\em not} 2-symmetric has always been, either implicitly or explicitly, made. The paper \cite{2sym} tries to 
fill in this gap
for the case of 2-symmetry and Lorentzian signature.

\section{Lorentzian 2-symmetry}
To deal with the problem of $k$-symmetric and 
$k$-recurrent spaces one needs to combine several different techniques. Among them (i) pure classical standard tensor 
calculus by using the Ricci and Bianchi identities; (ii) study of {\em parallel}
(also called covariantly constant) tensor and vector fields, 
and their implications on the manifold holonomy structure; and (iii) 
consequences on the curvature invariants. 
I now present the main points and results needed to reach the
sought results. It turns out that the 
so-called ``superenergy" and causal tensors \cite{S,BS} are 
very useful, providing positive quantities associated to tensors that can be used 
to replace the ordinary positive-definite metric available in proper 
Riemannian cases.

\subsection{Identities in 2-symmetric semi-Riemannian manifolds}
Of course, some tensor calculation is obviously needed, mainly to prove some helpful
quadratic identities. To start with, one needs a 
generalization of Proposition \ref{2implies1} to the case of non-generic points.
\begin{lemma}
Let $\espaitemps$ be an $n$-dimensional {\em 2-symmetric} 
semi-Riemannian manifold of any signature. If 
$\nabla_{\lambda}\nabla_{\mu}T_{\mu_1\dots \mu_q}=0$ then
\begin{eqnarray}
\sum_{i=1}^{q} \nabla_{\nu}R^{\rho}{}_{\alpha_i\lambda\mu}
T_{\alpha_1\dots\alpha_{i-1}\rho\alpha_{i+1}\dots\alpha_q}
-R^{\rho}{}_{\nu\lambda\mu}\nabla_{\rho}T_{\alpha_1\dots 
\alpha_q}=0,\label{eq4}\\
(\nabla_{\nu}R^{\rho}{}_{\tau\lambda\mu}+\nabla_{\tau}R^{\rho}{}_{\nu\lambda\mu})
\nabla_{\rho}T_{\mu_1\dots \mu_q}=0, \label{basic}\\
(\nabla_{\nu}R^{\rho}{}_{\mu}-\nabla_{\mu}R^{\rho}{}_{\nu})
\nabla_{\rho}T_{\mu_1\dots \mu_q}=0, \,\,\, 
(\nabla^{\rho}R_{\mu\nu}-2\nabla_{\nu}R^{\rho}{}_{\mu})
\nabla_{\rho}T_{\mu_1\dots \mu_q}=0. \label{treq4}
\end{eqnarray}
\end{lemma}
By using the decomposition of the Riemann tensor,
\begin{equation}
R_{\alpha\beta\lambda\mu}=C_{\alpha\beta\lambda\mu}+\frac{2}{n-2}\left(
R_{\alpha[\lambda}g_{\mu]\beta}-R_{\beta[\lambda}g_{\mu]\alpha}\right)-
\frac{R}{(n-1)(n-2)}\left(g_{\alpha\lambda}g_{\beta\mu}-
g_{\alpha\mu}g_{\beta\lambda}\right)\label{weyl}
\end{equation}
a selection of the formulas satisfied in 2-symmetric manifolds are 
given next
\begin{lemma}
The Riemann, Ricci and Weyl tensors of any $n$-dimensional {\em 
2-symmetric} semi-Riemannian manifold of any signature satisfy
\begin{eqnarray}
R^{\rho}{}_{\alpha\lambda\mu}R_{\rho\beta\gamma\delta}+
R^{\rho}{}_{\beta\lambda\mu}R_{\alpha\rho\gamma\delta}+
R^{\rho}{}_{\gamma\lambda\mu}R_{\alpha\beta\rho\delta}+
R^{\rho}{}_{\delta\lambda\mu}R_{\alpha\beta\gamma\rho}=0 \label{RicR}\\
R^{\rho}{}_{\nu\lambda\mu}\nabla_{\rho}R_{\alpha\beta\gamma\delta}+
R^{\rho}{}_{\alpha\lambda\mu}\nabla_{\nu}R_{\rho\beta\gamma\delta}+
R^{\rho}{}_{\beta\lambda\mu}\nabla_{\nu}R_{\alpha\rho\gamma\delta}+\nonumber\\
+R^{\rho}{}_{\gamma\lambda\mu}\nabla_{\nu}R_{\alpha\beta\rho\delta}+
R^{\rho}{}_{\delta\lambda\mu}\nabla_{\nu}R_{\alpha\beta\gamma\rho}=0 
\label{RicDR}\\
\nabla_{(\tau}R^{\rho}{}_{\nu)\lambda\mu}\nabla_{\rho}R_{\alpha\beta\gamma\delta}=0, 
\,\,
\nabla_{(\tau}R^{\rho}{}_{\nu)\lambda\mu}\nabla_{\rho}C_{\alpha\beta\gamma\delta}=0, 
\,\,
\nabla_{(\tau}R^{\rho}{}_{\nu)\lambda\mu}\nabla_{\rho}R_{\alpha\beta}=0, 
\label{basicR} \\
R_{\rho(\mu}R^{\rho}{}_{\nu)\alpha\beta}=0, \,\,\, 
R^{\rho}{}_{\mu[\alpha\beta}R_{\gamma]\rho}=0, \,\,\, 
C^{\rho}{}_{\mu[\alpha\beta}R_{\gamma]\rho}=0, \,\,\, 
R^{\rho\sigma}R_{\rho\mu\sigma\nu}=R_{\mu}{}^{\rho}R_{\rho\nu}, 
\label{RR}\\
R^{\rho}{}_{\alpha\lambda\mu}C_{\rho\beta\gamma\delta}+
R^{\rho}{}_{\beta\lambda\mu}C_{\alpha\rho\gamma\delta}+
R^{\rho}{}_{\gamma\lambda\mu}C_{\alpha\beta\rho\delta}+
R^{\rho}{}_{\delta\lambda\mu}C_{\alpha\beta\gamma\rho}=0, 
\label{RicC}\\
(n-2)\left(C_{\rho[\alpha}{}^{\lambda\mu}C^{\rho}{}_{\beta]\gamma\delta}+
C_{\rho[\gamma}{}^{\lambda\mu}C^{\rho}{}_{\delta]\alpha\beta}\right)-
2\left(R_{[\alpha}{}^{[\lambda}C^{\mu]}{}_{\beta]\gamma\delta}+
R_{[\gamma}{}^{[\lambda}C^{\mu]}{}_{\delta]\alpha\beta}\right)-\nonumber\\
-2\left(R_{\rho}{}^{[\lambda}\delta^{\mu]}_{[\alpha}C^{\rho}{}_{\beta]\gamma\delta}+
R_{\rho}{}^{[\lambda}\delta^{\mu]}_{[\gamma}C^{\rho}{}_{\delta]\alpha\beta}\right)+
2\frac{R}{n-1}\left(\delta^{[\lambda}_{[\alpha}C^{\mu]}{}_{\beta]\gamma\delta}+
\delta^{[\lambda}_{[\gamma}C^{\mu]}{}_{\delta]\alpha\beta}\right)=0 
\label{CC}
\end{eqnarray}
and their non-written traces, such as the appropriate specializations 
of (\ref{treq4}). Actually, (\ref{RicR}) and (\ref{RR}-\ref{CC}) are 
valid in arbitrary {\em semi-symmetric} spaces.
\end{lemma}

\subsection{Holonomy and reducibility in Lorentzian manifolds}
Some basic lemmas on local holonomy 
structure are also essential. The classical result here is the de Rham decomposition 
theorem \cite{deR,KN} for positive-definite metrics. However, this 
theorem does not hold as such for other signatures, and one has to 
introduce the so-called {\em non-degenerate reducibility}
\cite{Wu,Wu1,Wu2}. See also \cite{BI} for the particular case of Lorentzian 
signature. To fix ideas, recall that the holonomy group \cite{KN} 
of $\espaitemps$ is called reducible (when acting on the tangent 
spaces) if it leaves a non-trivial subspace of $T_p\varietat$ 
invariant. And it is called non-degenerately reducible if it leaves a 
non-degenerate subspace (that is, such that the restriction of the 
metric is non-degenerate) invariant.

Only a simple result is needed. This relates the 
existence of parallel tensor fields to the holonomy group of 
the manifold in the case of Lorentzian signature. It is a synthesis 
(adapted to our purposes) of the results in \cite{Hall} 
but generalized to arbitrary dimension $n$ (see \cite{2sym} for a proof): 
\begin{lemma}
\label{redornull}
Let $D\subset \varietat$ be a simply connected domain of an 
$n$-dimensional Lorentzian manifold $\espaitemps$ and assume that 
there exists a non-zero parallel symmetric tensor field 
$h_{\mu\nu}$ not proportional to the metric. Then $(D,g)$ is 
reducible, and further it is not non-degenerately reducible only if 
there exists a null parallel vector field which is the
unique parallel vector field (up to a constant of proportionality).
\end{lemma}
Some important remarks are in order here:
\begin{enumerate}
\item If there is a parallel 1-form $v_{\mu}$, then so is 
obviously $h_{\mu\nu}=v_{\mu}v_{\nu}$ and the manifold (arbitrary 
signature) is reducible, the Span of $v^{\mu}$ being invariant by the 
holonomy group. If $v_{\mu}$ is {\em not} null, then $\espaitemps$ is 
actually {\em non-degenerately} reducible. 
In this case, the metric can be decomposed into 
two orthogonal parts as $g_{\mu\nu}=c\, v_{\mu}v_{\nu} 
+(g_{\mu\nu}-c\, v_{\mu}v_{\nu})$, where $c=1/(v^{\mu}v_{\mu})$ is 
constant. Thus, necessarily 
$g_{\mu\nu}$ is a {\em flat extension} \cite{RWW} of a 
$(n-1)$-dimensional non-degenerate metric $g_{\mu\nu}-c\, 
v_{\mu}v_{\nu}$.
\item If there is a parallel non-symmetric tensor 
$H_{\mu\nu}$, then its symmetric part is also parallel, so that one can put $h_{\mu\nu}=H_{(\mu\nu)}$ in the 
lemma. In the case that $H_{\mu\nu}=H_{[\mu\nu]}\neq 0$ is 
antisymmetric, then in fact one can define 
$H_{\mu\rho}H_{\nu}{}^{\rho}=h_{\mu\nu}$, which is symmetric, 
parallel, non-zero and {\em not} proportional to the 
metric if $n>2$. For these last two statements, see e.g. \cite{BS}. 
\item Actually, the above can also be generalized to an arbitrary 
parallel $p$-form $\Sigma_{\mu_1\dots\mu_p}$ by defining 
$h_{\mu\nu}=\Sigma_{\mu\rho_2\dots\rho_p}\Sigma_{\nu}{}^{\rho_2\dots\rho_p}$.
\end{enumerate}

\subsection{Curvature invariants in 2-symmetric Lorentzian manifolds}
Recall that a curvature scalar invariant \cite{FKWC} is a 
scalar constructed polynomially from the Riemann tensor, the metric, 
the covariant derivative and possibly the volume element $n$-form of 
$\espaitemps$. They are called linear, quadratic, cubic, etcetera if 
they are linear, quadratic, cubic, and so on, on the Riemann tensor. 
This defines its {\em degree}. The {\em order} can be defined for 
homogeneous invariants, that is, so that they have the same number of 
covariant derivatives in all its terms. This number is the order of 
the scalar invariant. Of course, all non-homogeneous invariants can 
be broken into their respective homogeneous pieces, and therefore in 
what follows only the homogeneous ones will be considered. Similarly, 
one can define curvature 1-form invariants, or more generally, 
curvature $rank-r$ invariants in the same way but leaving 1, \dots , 
$r$ free indices \cite{FKWC}. 

A simple but very useful lemma is the following \cite{2sym}
\begin{lemma}
\label{nullorzero}
Let $(D,g)$ be as before with arbitrary signature. Any 
1-form curvature invariant which is parallel must be 
necessarily null (possibly zero).
\end{lemma}
It follows that, in 2-symmetric 
spaces, either $R$ is constant or $\nabla_{\mu} R$ is null and 
parallel. This is a particular example of the following 
general important result \cite{2sym}.
\begin{prop}
\label{constornull}
Let $D\subset \varietat$ be a simply connected domain of an 
$n$-dimensional {\em 2-symmetric} Lorentzian manifold $\espaitemps$. 
Then either
\begin{itemize}
\item all (homogeneous) scalar invariants of the Riemann tensor of 
order $m$ and degree up to $m+2$ are constant on $D$; or
\item there is a parallel null vector field on $D$.
\end{itemize}
\end{prop}
(Observe also that there will be no non-zero invariants 
involving derivatives of order higher than one. Then, the degree is 
necessarily greater or equal than the order.)

The previous proposition has immediate consequences providing more 
information about curvature invariants. For instance \cite{2sym}
\begin{coro}
\label{consequences}
Under the conditions of Proposition \ref{constornull}, either there 
is a parallel null vector field on $D$ or the following 
statements hold
\begin{enumerate}
\item All curvature scalar invariants of any order and degree formed 
as functions of the homogeneous ones of order $m$ and degree up to 
$m+2$ are constant on $D$;
\item All 1-form curvature invariants of order $m$ and degree up to 
$m+1$ are zero.
\item All scalar invariants with order equal to degree vanish.
\item All rank-2 tensor invariants with order equal to degree are 
zero.
\end{enumerate}
\end{coro}
{\bf Remark:} Of course, it can happen that the mentioned curvature 
invariants vanish {\em and} there is a parallel null 
vector field too.

There is a very long list of vanishing curvature invariants as a 
result of this Corollary---if there is no null parallel vector 
field---. The list of the quadratic ones is (only an 
independent set \cite{FKWC} is given,
omitting those contaning $\nabla_{\mu}R=0$):
\begin{eqnarray}
R^{\mu\nu}\nabla_{\alpha}R_{\mu\nu}=0,\,\,
R^{\mu\nu}\nabla_{\mu}R_{\nu\alpha}=0, \\
R^{\mu\nu\rho\alpha}\nabla_{\mu}R_{\nu\rho}=0,\,\,
R^{\mu\nu\rho\sigma}\nabla_{\mu}R_{\nu\rho\sigma\alpha}=0=
R^{\mu\nu\rho\sigma}\nabla_{\alpha}R_{\mu\nu\rho\sigma}, \label{r1}\\
\nabla_{\alpha}R^{\mu\nu}\nabla_{\beta}R_{\mu\nu}=
\nabla_{\mu}R_{\nu\beta}\nabla_{\alpha}R^{\mu\nu}=
\nabla_{\mu}R_{\nu\alpha}\nabla^{\mu}R^{\nu}{}_{\beta}=
\nabla_{\mu}R_{\nu\alpha}\nabla^{\nu}R^{\mu}{}_{\beta}=0, 
\label{r2ric}\\
\nabla^{\mu}R^{\nu\rho}\nabla_{\alpha}R_{\beta\rho\mu\nu}=
\nabla^{\mu}R^{\nu\rho}\nabla_{\mu}R_{\alpha\nu\beta\rho}=0,\\
\nabla_{\alpha}R^{\mu\nu\rho\sigma}\nabla_{\beta}R_{\mu\nu\rho\sigma}=
\nabla^{\sigma}R^{\mu\nu\rho\alpha}\nabla_{\sigma}R_{\mu\nu\rho\beta}=0\label{r2rie}
\end{eqnarray}
where of course the traces of (\ref{r2ric}-\ref{r2rie}) vanish, and 
one could also write the same expressions using the Weyl tensor 
instead of the Riemann tensor.

\section{Main results}
All necessary results to prove the main theorems have now been 
gathered. Then, by using the so-called future tensors and ``superenergy" techniques 
\cite{S,BS} one can prove the following\footnote{It must be stressed that this proof is only valid for Lorentzian manifolds, as the definition of future tensors requires 
this signature.} \cite{2sym}
\begin{theorem}
\label{ric-flat}
Let $D\subset \varietat$ be a simply connected domain of an 
$n$-dimensional {\em 2-symmetric} Lorentzian manifold $\espaitemps$. 
Then, if there is no null parallel vector field on $D$, 
$(D,g)$ is either Ricci-flat (i.e. $R_{\mu\nu}=0$) or locally symmetric.
\end{theorem}
Finally, one can at last prove that the narrow space left between locally 
symmetric and 2-symmetric Lorentzian manifolds can only be filled by 
spaces with a parallel null vector field. 
\begin{theorem}
Let $D\subset \varietat$ be a simply connected domain of an 
$n$-dimensional {\em 2-symmetric} Lorentzian manifold $\espaitemps$. 
Then, if there is no null parallel vector field on $D$, 
$(D,g)$ is in fact locally symmetric.
\end{theorem}
Thus we have arrived at
\begin{theorem}
\label{classi}
Let $D\subset \varietat$ be a simply connected domain of an 
$n$-dimensional {\em 2-symmetric} Lorentzian manifold $\espaitemps$. 
Then, the line element on $D$ is (possibly a flat extension of) the 
direct product of a certain number of locally symmetric proper 
Riemannian manifolds times either 
\begin{enumerate}
\item a Lorentzian locally symmetric 
spacetime (in which case the whole $(D,g)$ is locally symmetric), or 
\item a Lorentzian manifold with a parallel null vector 
field so that its metric tensor can be expressed locally as an 
appropriately restricted case of formula (\ref{brink}) below.
\end{enumerate}
\end{theorem}
Again, the following remarks are important:
\begin{enumerate}
\item Of course, the number of proper Riemannian symmetric manifolds 
can be zero, so that the whole 2-symmetric spacetime, if not locally 
symmetric, is given just by a line-element of the form (\ref{brink}) 
restricted to be 2-symmetric.
\item Although mentioned explicitly for the sake of clarity, it is 
obvious that the block added in any flat extension can also be 
considered as a particular case of a locally symmetric part building 
up the whole space.
\item This theorem provides a full characterization of the 
2-symmetric spaces using the classical results on the symmetric ones: 
their original classification (for the semisimple case) was given in 
\cite{Ber}, and the general problem was solved for Lorentzian 
signature in \cite{CW}. Combining these results with those for proper 
Riemannian metrics \cite{C,C1,H}, a complete classification is 
achieved.
\end{enumerate}

Thus, the only 2-symmetric non-symmetric Lorentzian manifolds contain 
a parallel null vector field. The most general 
local line-element for such a spacetime was discovered by Brinkmann 
\cite{Br} by studying the Einstein spaces which can be mapped 
conformally to each other. In appropriate local coordinates 
$\{x^0,x^1,x^i\}=\{u,v,x^i\}$, ($i,j,k,\dots =2,\dots ,n-1$) the 
line-element reads
\begin{equation}
ds^2=-2du(dv+Hdu+W_i dx^i)+g_{ij}dx^idx^j \label{brink}
\end{equation}
where the functions $H$, $W_i$ and $g_{ij}=g_{ji}$ are independent of 
$v$, otherwise arbitrary, and the parallel null vector field is given 
by
\begin{equation}
k_{\mu}dx^{\mu}=-du, \hspace{1cm} k^{\mu}\partial_{\mu}=\partial_v \, .
\end{equation}
It is now a 
simple matter of calculation to identify which manifolds among (\ref{brink})
are actually 2-symmetric. Using Theorem \ref{classi} and its remarks, 
this will provide ---by direct product with proper Riemannian symmetric 
manifolds if adequate--- all possible {\em non-symmetric} 2-symmetric 
spacetimes. By doing so \cite{2sym} one finds, among other results, that (i) the $g_{ij}$ are a one-parameter family, depending on $u$, of locally symmetric proper Riemannian metrics\footnote{As these are classified 
in e.g. \cite{C1,H}, the part $g_{ij}$ of the metric is completely 
determined. For an explicit formula, one only has to take any of them 
from the list and let any arbitrary constants appearing there to be 
functions of $u$.}; (ii) for a given choice of $g_{ij}$ in agreement with the previous point, the integrability conditions provide the explicit form of the functions $H$ and $W_i$; (iii) finally, the scalar curvature coincides with the corresponding scalar 
curvature $\bar{R}$ of $g_{ij}$: $R=\bar{R}$.
Due to (i), the function $\bar{R}$ depends only on $u$, and thus 
the 2-symmetry implies
\begin{equation}
R=\bar{R}(u)=a\, u + b \label {scalar}
\end{equation}
where $a$ and $b$ are constants. In particular, 
$\nabla_{\mu}R=-ak_{\mu}$. Thus, we see that given any locally symmetric proper Riemannian $g_{ij}$ and letting the constants appearing there to be functions of $u$ is too general, and these functions are restricted by the 2-symmetry so that, for example, (\ref{scalar}) holds.


\begin{thebibliography}{99}
\bibitem{BI} L. B\'erard-Bergery and A. Ikemakhen, ``On the holonomy of 
Lorentzian manifolds", in
{\it Differential geometry: geometry in mathematical physics and 
related topics, Los Angeles, CA, 1990} pp. 27--40, {\it Proc. Sympos. 
Pure Math.}, {\bf 54}, Part 2, (Amer. Math. Soc., Providence, R.I., 
1993)
\bibitem{Ber} M. Berger, {\it Ann. Sci. \'{E}cole Norm. Sup.} 
{\bf 74} 85--177 (1957)
\bibitem{BS} G. Bergqvist and J.M.M. Senovilla,  {\it Class. Quantum Grav.} 
{\bf 18} 5299-5325 (2001)
\bibitem{BKV} E. Boeckx, O. Kowalski, L. Vanhecke, {\it Riemannian 
manifolds of conullity two} World Sci. Singapore 1996
\bibitem{Br} H. W. Brinkmann, {\it Math. Ann.} {\bf 94} 119-145 (1925)
\bibitem{CM} M. Cahen and R. McLenaghan,  {\em C. R. Acad. 
Sci. Paris S\'er. A-B} {\bf 266} A1125--A1128 (1968)
\bibitem{CP} M. Cahen and M. Parker,  {\em Bull. Soc. Math. Belg.}  {\bf 
22} 339--354 (1970)
\bibitem{CP1} M. Cahen and M. Parker,  {\it Mem. Amer. Math. Soc.} {\bf 24} no. 229 (1980)
\bibitem{CW} M. Cahen and N. Wallach,
{\em Bull. Amer. Math. Soc.} {\bf 76} 585--591 (1970)
\bibitem{C} \'{E}. Cartan, {\em Bull. Soc. Math. France} {\bf 54} 214-264 (1926); {\bf 
55} 114-134 (1927)
\bibitem{C1} \'{E}. Cartan, {\it Le\c{c}ons sur la G\'eom\'etrie des 
Espaces de Riemann}, 2nd ed., Gauthier-Villars, Paris 1946
\bibitem{CS1} C.D. Collinson and F. S\"oler, {\it Tensor (N.S.)}  {\bf 30} 
87--88 (1976)
\bibitem{FKWC} S.A. Fulling, R.C.  King, B.G. Wybourne, C.J. 
Cummins,
{\it Class. Quantum Grav.} {\bf 9} 1151--1197 (1992)
\bibitem{Hall} G.S. Hall, {\it J. Math. Phys.} {\bf 32} 
181--187 (1991)
\bibitem{H} S. Helgason, {\it Differential Geometry, Lie groups, and 
Symmetric Spaces}, Academic Press, New York 1978
\bibitem{KN} S. Kobayashi and K. Nomizu, {\it Foundations of Differential 
Geometry}, Wiley, Interscience, New York, vol.I 1963, vol.II 1969
\bibitem{L} A. Lichnerowicz,  Courbure, nombres de Betti, et espaces 
sym\'etriques, {\it Proceedings of the International Congress of 
Mathematicians, Cambridge, Mass., 1950} vol. 2,  216--223
(Amer. Math. Soc., Providence, R. I., 1952)
\bibitem{NO} K. Nomizu and H. Ozeki, {\it Proc. Nat. Acad. Sci. USA} {\bf 48} 206-207 (1962)
\bibitem{O} B. O'Neill, {\it Semi-Riemannian Geometry}, Academic Press, New 
York 1983
\bibitem{deR} G. de Rham, {\it Comment. Math. Helv.} {\bf 26} 328--344 (1952)
\bibitem{RWW} {\it H.S. Ruse, A.G. Walker, T.J. Willmore, Harmonic spaces}, 
(Consiglio Nazionale delle Ricerche Monografie Matematiche, 8)
Edizioni Cremonese, Rome 1961
\bibitem{S} J.M.M. Senovilla, {\it Class. 
Quantum Grav.} {\bf 17} 2799-2841 (2000)
\bibitem{2sym} J.~M.~M.~Senovilla, preprint (2005) (available upon request).
\bibitem{Sz} Z.I. Szab\'o, {\it J. Diff. Geom.} {\bf 17} 531-582 (1982)
\bibitem{Sz1} Z.I. Szab\'o, {\it Geom. Dedicata} {\bf 19} 65-108 (1985)
\bibitem{Tak} H. Takeno, {\it Tensor (N.S.)} {\bf  27} 309--318  (1973)
\bibitem{Ta} S. Tanno, 
{\it  Ann. Mat. Pura Appl. (4)}  {\bf 96} 233--241 (1972)
\bibitem{Wu} H. Wu, {\it Illinois J. Math.} {\bf 8} 291--311 (1964)
\bibitem{Wu1} H. Wu,  {\it Bull. Amer. Math. Soc.} {\bf 70} 610--617 (1964)
\bibitem{Wu2} H. Wu, {\it Pacific J. Math.} {\bf 20} 351--392 (1967)

\end{thebibliography}
\end{document}